\documentclass[12pt]{article}   
\usepackage{amsmath, amssymb}  
\usepackage{amsthm}

\usepackage{verbatim}
\usepackage{amsfonts}
\usepackage{amssymb}
\usepackage{mathrsfs}
\usepackage[dvips]{graphicx}
\usepackage{pstricks} 
\usepackage{pst-node}

\usepackage{setspace}
\usepackage[inner=2.4cm,outer=4cm,top=2.4cm,bottom=2.4cm,marginparwidth=90pt]{geometry}   
\usepackage{times}

\usepackage[notref,notcite,final]{showkeys}  
  
\numberwithin{equation}{section}

\usepackage{framed}  

\newcommand{\bigR}{{\mathbb R}}

\newtheorem{theorem}{Theorem}[section]  
\newtheorem{lemma}[theorem]{Lemma}  
\newtheorem{proposition}[theorem]{Proposition}  
  
\newtheorem{remark}[theorem]{Remark}

\allowdisplaybreaks[1]  

 \title{On the extension of axially symmetric volume preserving mean curvature flow}
\author{John Head, Sevvandi Kandanaarachchi} 

\date{}
  
\begin{document}  
  
\maketitle

 \onehalfspacing  
  

\abstract
We study the provenance of singularity formation under mean curvature flow and volume preserving mean curvature flow in an axially symmetric setting.  We prove that if the mean curvature is uniformly bounded on any finite time interval, then no singularities can develop during that time under both mean curvature flow and volume preserving mean curvature flow.

\section{Introduction}
\noindent
Consider a smooth, $n$-dimensional hypersurface immersion $\mathbf{x}_0: M^n \rightarrow \mathbb{R}^{n+1}$. The solution of mean curvature flow generated by $\mathbf{x}_0(M^n)$ is the one-parameter family $\mathbf{x}: M^n \times [0,T) \rightarrow \mathbb{R}^{n+1}$ of smooth immersions satisfying
\begin{equation}\label{eq:int_2}
\frac{\partial}{\partial t} \mathbf{x}(l, t)=-H(l,t)\nu(l,t), \hspace{5 mm} l\in M^n, t\geq0 \, ,
\end{equation}
\noindent
with $\mathbf{x}(\cdot,0) = \mathbf{x}_0$. Here $\nu(l,t)$ represents a choice of unit normal (the outward-pointing unit normal in the closed setting), and $H(l,t)$ is the mean curvature. According to our choice of signs, the right-hand side is the mean curvature vector and the mean curvature of the round sphere is positive. We henceforth write $M_t = \mathbf{x}(\cdot, t)(M^n)$.
 \\

\noindent
If the evolving surfaces $M_t$ are assumed to be compact and enclose a prescribed volume $V$, this process is governed by the evolution equation
\begin{equation}\label{eq:int_1}
\frac{d}{dt} \mathbf{x}(l,t) = -\left( H(l,t) -h(t)\right)\nu(l,t), 
\hspace{5 mm} l\in M^n, t\geq0,
\end{equation}
\noindent
where
\[h(t) =\frac{\int_{M_t} H d\mu_t}{\int_{M_t} d\mu_t} \]
is the average of the mean curvature and $d\mu_t$ denotes the surface measure on $M_t$.
\\

\noindent
As initial data we select a smooth, compact, 2-dimensional hypersurface $M_0$ with boundary $\partial M_0 \neq \emptyset$. We assume that $M_0$ is smoothly embedded in the domain
$$ G = \{ \mathbf{x} \in \bigR^{3} : a < x_1 < b \} \, , \hspace{3mm}  0<a<b \, ,
$$
and that the free boundary satisfies $\partial M_0 \subset \partial G$. Moreover, we consider an axially symmetric surface contained in the region $G$ between the two parallel planes $x_1 = a$ and $x_1 = b$.\\

\noindent
Motivated by the fact that the stationary solution to the associated Euler Lagrange problem satisfies a Neumann boundary condition, we assume that the surface meets the planes $x_1 = a$ and $x_1 = b$ orthogonally. In this setting, we consider the general question of whether a singularity can develop if the mean curvature remains bounded. \\

\noindent
Le and Sesum studied closely related problems for mean curvature flow in \cite{LeSe10}. In particular, they proved that if all singularities are of type I, the mean curvature blows up at the first singular time $T$ . The major difference between this result and what is discussed here is that we work with axially symmetric surfaces but do not have any restrictions on the type of singularity.  Li and Wang study related problems in  \cite{li2016extension}.  \\

\noindent
This paper is organised as follows. In sections \ref{Section_Evolution_Eqs_Vol} and  \ref{Section_VolFlow},  we study volume preserving mean curvature flow and investigate whether a singularity can develop if $H$ is bounded. In sections \ref{Section_Evolution_Eqs_MCF} and \ref{Section_MCF} we study the corresponding question for mean curvature flow. We prove the following two theorems.

\begin{theorem} \label{Theorem_HBdd}
Consider a smooth solution $M_t$ of \eqref{eq:int_1} on $[0,T)$ for some $T>0$. Suppose that $ |H(l,t)| \leq c_0 \,$ for all $l \in M^2$  and for all $t \in [0, T)$. Then there exists a constant $C>0$ depending only on $c_0$, $M_0$ and $T$ such that $|A|(l,t) \leq C$ for all $t \in [0, T).$ That is, the flow can be extended past time $T$.
\end{theorem}

\begin{theorem} \label{Theorem_HBdd2}
Consider a smooth solution $M_t$ of \eqref{eq:int_2} on $[0,T)$ for some $T>0$. Suppose that $ |H(l,t)| \leq c_0' \,$ for all $l \in M^2$  and for all $t \in [0, T)$. Then there exists a constant $C'>0$ depending only on $c_0'$ and $M_0$ such that $|A|(l,t) \leq C'$ for all $t \in [0, T).$ That is, the flow can be extended past time $T$.
\end{theorem}

\noindent

\noindent 
{\small{{\bf Acknowledgments.} This work forms part of the second author's PhD dissertation. She would like to thank her advisor Maria Athanassenas for her guidance and support throughout the preparation of this work. The authors would also like to thank Gerhard Huisken for valuable discussions on this work.}}

\section{Notation}\label{Section_Notation}

\noindent
We follow the notation used by the second author in \cite{SKThesis}, which agrees with the notation adopted by Huisken in \cite{GH2} and by Athanassenas in \cite{MA1}. Let $\rho_0:[a,b] \rightarrow  \mathbb{R} $ be a smooth, positive function on the bounded interval $[a,b]$ with $\rho'_0(a) = \rho'_0(b)=0 $. Consider the $2$-dimensional hypersurface $M_0$ in $\mathbb{R}^{3} $ generated by rotating the graph of $\rho_0$ about the $x_1$-axis. We evolve $M_0$ along its mean curvature vector while preserving the enclosed volume subject to Neumann boundary conditions at $x_1 = a$ and $x_1 = b$. Equivalently, we can consider the evolution of a periodic surface defined on the entire $x_1$ axis. By definition the evolution preserves axial symmetry.\\

\noindent
The position vector $\mathbf{x}$ satisfies
\begin{equation}\label{notEq:1.3}
\frac{\partial}{\partial t}\mathbf{x} = -(H-h)\nu = \mathbf{H} +h\nu =  \Delta{\mathbf{x}}+h\nu \, ,
\end{equation}
\noindent 
where $\mathbf{H}$ denotes the mean curvature vector. \\

\noindent
We write $\mathbf{i}_1, \mathbf{i}_2, \mathbf{i}_{3}$ for the standard basis vectors in $\mathbb{R}^{3}$ corresponding to the $x_1, x_2, x_{3}$ axes respectively. We then introduce a local orthonormal frame  $\tau_1(t), \tau_2(t) $ on $M_t$ such that 
\[ \left\langle\tau_2(t), \mathbf{i}_1\right\rangle = 0,  \hspace{5mm} \text{and} \hspace{5mm} \left\langle\tau_1(t), \mathbf{i}_1\right\rangle > 0  \, .\]

\noindent
We denote by $ \omega= \frac{\hat{\mathbf{x}}}{|\hat{\mathbf{x}}|} \in \mathbb{R}^{3}$ the outward-pointing unit normal to the cylinder intersecting $M_t$ at the point $\mathbf{x}(l,t)$. Here $\hat{\mathbf{x}}= \mathbf{x} - \left\langle \mathbf{x}, \mathbf{i}_1 \right\rangle \mathbf{i}_1$.  Using this, we define 

$$y = \left\langle \mathbf{x}, \omega \right\rangle \hspace{3 mm} \text{and}  \hspace{3mm} v=\left\langle \omega, \nu \right\rangle ^{-1}\, .$$ 
We henceforth refer to $y$ as the { \it height function} and $v$ as the {\it gradient function}. We have $\rho: [a, b] \times [0, T) \rightarrow \bigR$, whereas $y: M^2 \times [0, T) \rightarrow \bigR $. Note that $v$ is a geometric quantity related to the inclination angle. More precisely, $v$ corresponds to $\sqrt{1+ \rho'^2}$ in the axially symmetric setting. The study of this quantity has given rise to fundamental gradient estimates in the graphical setting, see for example \cite{EH89,Eck97Mink}. \\

\noindent
We denote by $g=\{g_{ij}\}$ the induced metric and by $A=\{h_{ij}\}$ the second fundamental form at the space-time point $(l,t)\in M^2\times [0,T)$. \\

\noindent
Finally, we introduce the quantities (see also \cite{GH2})
\begin{equation}\label{notEq:1.1}
 p =  \left\langle\tau_1, \mathbf{i}_1 \right\rangle y^{-1},  \hspace{10 mm} q= \left\langle\nu, \mathbf{i}_1 \right\rangle y^{-1}, 
\end{equation}

\noindent
which satisfy
\begin{equation}\label{notEq:1.2}
p^2 + q^2 = y^{-2} \, .
\end{equation}

\noindent
The second fundamental form has one eigenvalue given by
$$p = \frac{1}{\rho \sqrt{1+ \rho'^2}}$$
and the other equal to
\[ k= \left\langle \overline{\nabla}_{1} \nu, \tau_1 \right\rangle = \frac{-\rho''}{(1+\rho'^2)^{3/2}}. \]

\section{Preliminary results - Volume Preserving MCF }\label{Section_Evolution_Eqs_Vol}
 
\noindent

In this section we prove some preliminary estimates which will form the foundation of the rescaling and contradiction argument which follows in the next section. We begin by recalling the well-known evolution equations.

\begin{lemma}\label{Lemma_Evolution_Equations} We have the following evolution equations: 
\begin{itemize}
\item[(i)] $\frac{\partial}{\partial t} y = \Delta y - \frac{1}{y} + hpy \, ; $
\item[(ii)]$ \frac{\partial}{\partial t} v = \Delta v  -|A|^2v + \frac{v}{y^2} - \frac{2}{v}|\nabla v|^2\, ;  $
\item [(iii)] $ \frac{\partial}{\partial t} k = \Delta k + |A|^2k -2q^2(k-p) - hk^2 \, ;  $
\item [(iv)] $ \frac{\partial}{\partial t} p = \Delta p + |A|^2p + 2q^2(k-p) - hp^2 \, ;  $
\item [(v)]$ \frac{\partial}{\partial t} H = \Delta H + (H -h) |A|^2 \, .$
\end{itemize}

\end{lemma}
\begin{proof}
\noindent
Evolution equations (i), (ii), (iii) and (iv) are proved in \cite{AthKan1}. Equation (v) is proved in \cite{GH3}. 
\end{proof}

\noindent
We seek to proceed using the maximum principle. Estimates on $h$ are therefore essential as this is a global term. In the axially symmetric setting with Neumann boundary data, these were established by Athanassenas in \cite{MA2}:

\begin{lemma} \label{Maria's Prop 1.4} Consider a smooth solution $M_t$ of \eqref{eq:int_1} on $t \in [0, T)$ for some $T>0$. Then there exist constants $c_1,c_2>0$ depending only on the initial hypersurface $M_0$ such that the mean value $h$ of the mean curvature satisfies
$$0 < c_1 \leq  h \leq c_2 \, . $$
\end{lemma}

\noindent
We use this to produce a height dependent gradient estimate. 

\begin{lemma} \label{Lemma_UV} Consider a smooth solution $M_t$ of \eqref{eq:int_1} on $t \in [0, T)$ for some $T>0$. Then there exists a constant $c_3>0$ depending only on $M_0$ and $T$ such that
$$ vy \leq c_3 $$
for all $t \in [0, T)$.
\end{lemma}

\begin{proof}
From Lemma \ref{Lemma_Evolution_Equations}, for any $c>0$ we have
\[ \frac{\partial}{\partial t} (yv -c t) = \Delta(yv) -\frac{2}{v}\left\langle \nabla v \, , \nabla (yv) \right\rangle - yv|A|^2 + h-c  .\]
Applying Lemma \ref{Maria's Prop 1.4} we obtain
\begin{align*}
yv-c_2 t &\leq \max_{M_0} yv \, .
\end{align*}
This establishes the claim.
\end{proof}

\begin{lemma} \label{Prop_KonP}
Consider a smooth solution $M_t$ of \eqref{eq:int_1} on $t \in [0, T)$ for some $T>0$. Then there exists a constant $c_4>0$ depending only on $M_0$ such that the principal curvatures $k$ and $p$ satisfy
$$\frac{k}{p} < c_4$$
for all $t\in[0,T)$.
\end{lemma}

\begin{proof}
Using Lemma \ref{Lemma_Evolution_Equations} we calculate
\[ \frac{\partial}{\partial t} \left( \frac{k}{p} \right) = \Delta \left( \frac{k}{p} \right)+
	\frac{2}{p} \left\langle \nabla p \, , \nabla \left(\frac{k}{p}\right) \right\rangle +
	2\left(\frac{q}{p}\right)^2\left(p^2 - k^2\right) +
	\frac{hk}{p}\left(p-k\right) \, .\]
See also equation $(19)$ in \cite{GH2}.
If $\frac{k}{p} \geq 1$ then $p-k\leq 0$ and from the parabolic maximum principle we obtain
\begin{equation}
\frac{k}{p} \leq \max \left\{1, \max_{M_0} \frac{k}{p} \right\}  \, .
\end{equation}
This completes the proof.
\end{proof}

\noindent
Finally, we establish an estimate for $\frac{|k|}{p}$ assuming a uniform bound on the mean curvature.

\begin{proposition} \label{Prop_KonP3}
Consider a smooth solution $M_t$ of \eqref{eq:int_1} on $t \in [0, T)$ for some $T>0$.
Suppose that there exists a constant $c_0>0$ such that $|H(l,t)| \leq c_0$ for all $l\in M^2$ and for all $t\in[0,T)$.
Then there exists a constant $c_5>0$ depending only on $c_0$, $M_0$ and $T$ such that
$$\frac{|k|}{p} \leq c_5$$
for all $t \in[0, T)$. 
\end{proposition}

\begin{proof} Since
$$ -c_0 \leq H = k + p \leq c_0  $$
and $\frac{1}{p} = vy$, we obtain from Lemma \ref{Lemma_UV} that
$$ \frac{|k|}{p} \leq  1+ c_0c_3 \, ,$$
as required.
\end{proof}

\noindent
In the next section we combine these estimates with standard rescaling techniques to prove Theorem 1.1.

\section{Rescaling Volume Preserving MCF} \label{Section_VolFlow}

We begin by introducing a rescaling procedure similar to those which appear frequently in the literature, see for example \cite{HS2}.
Given the smooth solution $M_t$ of \eqref{eq:int_1} for $t \in [0, T)$, we consider integers $i\geq 1$, times $t_i\in[0 , T - \frac{1}{i} ]\, ,$ and points $l_i \in M^2$ such that the following hold: $t_i\rightarrow T \,$, $\mathbf{x}(l_i, t_i)$ lies on the $\left(x_1, x_{3}\right)$ plane and
\begin{equation}\label{eq_Hbdd_eq15}
  |A|(l_i, t_i)  = \max_{\substack{ l \in M^{2} \\ t \leq T- \frac{1}{i} } } |A|(l,t) \,    .
\end{equation}
Let $\alpha_i = |A|(l_i, t_i) $ and $ \mathbf{x}_i = \mathbf{x}(l_i, t_i) \, .$ We now perform a parabolic rescaling of $M_t$, giving rise to the family $\tilde{M}_{i, \tau}$ defined by
\begin{equation} \label{HbddEq_1}
 \tilde{\mathbf{x}}_i(l \,, \tau) = \alpha_i \left( \mathbf{x}(l \,, \alpha_i^{-2} \tau + t_i )	- \left\langle \mathbf{x}_i , \mathbf{i}_1 \right\rangle \mathbf{i}_1 \right)  \, . 
\end{equation} 
Observe that $ \tau \in [ -\alpha_i^2 t_i, \alpha_i^2 (T- t_i - \frac{1}{i} )] \, .$ \\

\noindent
By construction we rescale from a point on the axis of rotation corresponding to the maximum curvature, preserving axial symmetry. We denote by $\tilde{\rho}_{i , \tau}$ the generating curves of $\tilde{M}_{i, \tau} \,$ (the cross section of $\tilde{M}_{i, \tau}$ in the $\left(x_1, x_3\right)$ plane). In addition, we denote by $|\tilde{A}_i| $ and $\tilde{H_i}$ the second fundamental form and mean curvature of $\tilde{M}_{i, \tau} \,$, respectively. By definition
$$  \tilde{H_i}( \cdot \,, \tau ) = \alpha_i^{-1} H(\cdot \,, \alpha_i^{-2} \tau + t_i ) \,  \hspace{3mm} \text{and} \hspace{3mm} |\tilde{A}_i|( \cdot \,, \tau ) = \alpha_i^{-1}	|A| (\cdot \,, \alpha_i^{-2} \tau + t_i )  \, .$$					
For $t \leq T - \frac{1}{i}$ we have 
\begin{equation}\label{HbddEq_2}
 \alpha_i^{-1}	|A| (\cdot \,, \alpha_i^{-2} \tau + t_i ) \leq 1 \, . 
\end{equation}

\noindent
In order to verify that each rescaled flow preserves volume, note that
$$  d\tilde{\mu}_\tau = \alpha_i^2 d\mu_t  \, ,$$
where $d\tilde{\mu}_\tau$ denotes the surface measure on $\tilde{M}_{i, \tau}$. Now let
\begin{equation*}
 \tilde{h}_i(\tau) = \frac{ \int_{\tilde{M}_{i, \tau}} \tilde{H}_i(\cdot \, , \tau) d\tilde{\mu}_{\tau}}{ \int_{\tilde{M}_{i, \tau}} d\tilde{\mu}_{\tau}} \, ,
\end{equation*}
\noindent
from which we obtain
$$ \tilde{h}_i (\tau) = \frac{ \int_{M_t} \alpha_i^{-1}  H (\cdot \,, \alpha_i^{-2}\tau + t_i) \alpha_i^2 d\mu_t}{ \int_{M_t} \alpha_i^2 d\mu_t} \, = \alpha_i^{-1} h(\alpha_i^{-2}\tau + t_i) = \alpha_i^{-1} h(t) \,  , $$

\noindent
and
\begin{equation}\label{HbddEq_3}
\frac{\partial}{\partial \tau} \tilde{\mathbf{x}}_i 		= - \alpha_i^{-1} (H- h ) \nu		= - (\tilde{H}_i - \tilde{h}_i )\nu \, .
\end{equation}

\noindent We now prove Theorem \ref{Theorem_HBdd}. \\

\begin{proof} [Proof of Theorem \ref{Theorem_HBdd}]
Suppose that $|A|^2 \rightarrow \infty$ as $t \rightarrow T$, i.e. a singularity forms at time $t=T$. We rescale the flow using $\eqref{HbddEq_1}$. As a first step, we show that the rescalings cannot drift away to infinity. Applying Proposition \ref{Prop_KonP3} we can find a constant $c_6$ depending only on $c_0$, $M_0$ and $T$ such that
$$ |A|=\sqrt{k^2+p^2}\leq c_6 p \leq c_6 y^{-1} . $$
After rescaling, this becomes
$$ |\tilde{A}_i| \leq c_6 (\alpha_{i} y)^{-1} =c_6 \tilde{y}_i^{-1} . $$
Since $|\tilde{A}_i|(l_i,0)=1$ for all $i$, we therefore have a bound on $\tilde{y}$ and can extract a convergent subsequence of points on the $x_3$ axis.
\\

\noindent
Along the sequence of rescaled flows we have the uniform curvature bound $|\tilde{A}_i|^2 \leq 1 $. Since each rescaled flow again satisfies (\ref{HbddEq_3}), this gives rise to uniform bounds on all covariant derivatives of the second fundamental form, see for example \cite{GH3}. By standard methods, based on the Arzela-Ascoli theorem, we can therefore find a further subsequence which converges uniformly in $C^{\infty}$ on compact subsets of $\mathbb{R}^{3}\times\mathbb{R}$ to a non-empty smooth limit flow which exists on an interval $(-\infty,\beta)$ where $\beta\in[0,\infty]$. We now analyse the properties of this limit flow, which we label $\tilde{M}_{\infty,\tau}$. \\

\noindent
Since $|H|$ is bounded by assumption and $h$ is bounded by Lemma \ref{Maria's Prop 1.4}, we obtain
$$\lim_{i\to\infty}\tilde{H}_i=0 \quad \mbox{and} \lim_{i\to\infty}\tilde{h}_i=0 \, .$$
Thus $\tilde{M}_{\infty,\tau}$ is a stationary solution of the flow, which we relabel $\hat{M}$, and must therefore be the catenoid. \\

\noindent
We henceforth use a 'hat' to indicate that a geometric quantity is associated with the catenoid.  The catenoid is obtained by rotating $\hat{y} = c_7 \cosh(c_7^{-1}\hat{x}_1) $ around the $x_1$ axis. For any $\epsilon>0$ and for any $l \in M^2$ we can find $I_0\in \mathbb{N}$ such that for any fixed $\tau_0\in (-\alpha_{I_0}^2t_{I_0},0)$ we have
$$ \hat{v}(l)\hat{y}(l) -\epsilon  \leq  \tilde{v}_i(l, \tau_0) \tilde{y}_i(l, \tau_0)  \hspace{5mm} \text{for all} \hspace{2mm} i >I_0 \, . $$

\noindent
On the catenoid, $\hat{v} = \sqrt{1 + \hat{y}'^2}= \cosh(c_7^{-1}\hat{x}_1)$. It therefore follows from Lemma \ref{Lemma_UV} that
\begin{equation*}
 \frac{c_7}{ 2\alpha_i}\left( \cosh (2c_7^{-1}\hat{x}_1) +1  \right) - \frac{\epsilon}{ \alpha_i} \leq c_3  \hspace{5mm} \text{for all} \hspace{2mm} i >I_0 \, .
\end{equation*}

\noindent
For fixed $i$, the left-hand side can be made as large as we like, giving rise to a contradiction. \\

\noindent
We can therefore find a constant $C$ such that $|A| \leq C$ for all $t \in [0, T)$. Again applying Theorem $4.1$ in \cite{GH3}, we obtain estimates on all covariant derivatives of $|A|$, allowing us to extend the flow beyond $T$. This completes the proof.  

\end{proof}

\section{Preliminary results - MCF }\label{Section_Evolution_Eqs_MCF}

\noindent
As before, we consider an axially symmetric surface with Neumann boundary data in $\bigR^{3}$. We follow the notation used in the volume preserving case, but study the equation (\ref{eq:int_2}).
Without imposing any additional restrictions, we establish that a singularity cannot develop if $|H|$ remains bounded. \\

\noindent
It is well known that if $H>0$ then $|A|^2 \leq c H^2 $ and therefore no singularity can develop, see \cite{GH2}. However, this argument is based on the evolution equation for $\frac{|A|^2}{H^2}$ and does not apply if the mean curvature changes sign. Since we do not assume that the mean curvature is positive on $M_0$, we instead use the blow up techniques employed in the volume-preserving setting in the previous section.  \\

\noindent
We begin by recalling the relevant evolution equations.

\begin{lemma}\label{Lemma_Evolution_Equations_2} We have the following evolution equations: 
\begin{itemize}
\item[(i)] $\frac{\partial}{\partial t} y = \Delta y - \frac{1}{y}  \, ; $
\item[(ii)]$ \frac{\partial}{\partial t} v = \Delta v  -|A|^2v + \frac{v}{y^2} - \frac{2}{v}|\nabla v|^2\, ;  $
\item [(iii)] $ \frac{\partial}{\partial t} k = \Delta k + |A|^2k -2q^2(k-p)  \, ;  $
\item [(iv)] $ \frac{\partial}{\partial t} p = \Delta p + |A|^2p + 2q^2(k-p)  \, ;  $
\item [(v)]$ \frac{\partial}{\partial t} H = \Delta H + H  |A|^2 \, . $
\end{itemize}

\end{lemma}
\begin{proof}
\noindent
Evolution equations (i), (iii) and (iv)  are proved in \cite{GH2}, (ii) appears in \cite{JBT1} and (v) can be found in  \cite{GH1}.  
\end{proof}

\noindent
The next result is an analogue of Lemma \ref{Lemma_UV} for mean curvature flow.
 
\begin{lemma} \label{Lemma_UV2} Consider a smooth solution $M_t$ of \eqref{eq:int_2} on $[0,T)$ for some $T>0$. Then there exists a constant $c_8>0$ depending only on the initial hypersurface $M_0$ such that
$$ vy\leq c_8 $$
for all $t\in[0,T)$.
\end{lemma}

\begin{proof}
From Lemma \ref{Lemma_Evolution_Equations_2} we compute
\[ \frac{\partial}{\partial t} (yv)  = \Delta(yv) -\frac{2}{v}\left\langle \nabla v \hspace{1mm}, \nabla (yv) \right\rangle - yv|A|^2  \, .\]
The maximum principle yields
$$yv  \leq \max_{M_0} yv \, . $$
\end{proof}

\noindent
We now record analogues of Lemma \ref{Prop_KonP} and Proposition \ref{Prop_KonP3} in our current setting.

\begin{lemma} \label{MC_Prop_KonP2} Consider a smooth solution $M_t$ of \eqref{eq:int_2} on $[0,T)$ for some $T>0$. Then there exists a constant $c_9>0$ depending only on $M_0$ such that
$$ \frac{k}{p}\leq c_9 $$
for all $t\in[0,T)$.
\end{lemma}

\begin{proof}
This is proved in \cite{GH2}.

\end{proof}

\begin{proposition} \label{MC_Prop_KonP2} Consider a smooth solution $M_t$ of \eqref{eq:int_2} on $[0,T)$ for some $T>0$. Suppose that there exists a constant $c_0'>0$ such that $|H(l,t)|\leq c_0'$ for all $l\in M^2$ and for all $t\in[0,T)$. Then there exists a constant $c_{10}>0$ depending only on $c_0'$ and $M_0$ such that
$$ \frac{|k|}{p}\leq c_{10} $$
for all $t\in[0,T)$.
\end{proposition}

\begin{proof}
Using Lemma \ref{Lemma_UV2}, the argument proceeds in the same way as the proof of Proposition \ref{Prop_KonP3}.
\end{proof}

\noindent
These estimates can now be used to prove Theorem \ref{Theorem_HBdd2}.

\section{Rescaling MCF}\label{Section_MCF}
\noindent
We now have all of the prerequisites to prove Theorem \ref{Theorem_HBdd2}. We will employ the rescaling procedure defined in \eqref{HbddEq_1}. Note, however, that
\begin{equation*}
\frac{\partial}{\partial \tau} \tilde{\mathbf{x}}_i  		= - \alpha_i^{-1} H \nu \, 
										= - \tilde{H}_i \nu \, ,
\end{equation*}
giving rise to a sequence of axially symmetric rescaled mean curvature flows.

\noindent

\begin{proof} [Proof of Theorem \ref{Theorem_HBdd2}]
This follows from an argument which is identical to the proof of Theorem \ref{Theorem_HBdd}. We emphasize that in the current setting we move from the uniform curvature bound $|\tilde{A}_i|^2\leq 1$ to uniform bounds on all covariant derivatives of $\tilde{A}_i$ via Proposition $2.3$ in \cite{GH2}.   
\end{proof}

\bibliographystyle{acm}
\bibliography{citations}

\vspace{.5cm}
\noindent
School of Mathematical Sciences, \\
Monash University, Victoria, Australia  \\
john.head@monash.edu \\

\noindent
Department of Econometrics and Business Statistics,   \\
Monash University, Victoria, Australia \\
sevvandi.kandanaarachchi@monash.edu

\end{document}